# Use of Homotopy Perturbation Method for Solving Multi-point Boundary Value Problems

[1, a] Shahid S. Siddiqi and [1, 2, b] Muzammal Iftikhar
[1]Department of Mathematics, University of the Punjab, Lahore 54590, Pakistan
[2]Department of Mathematics, University of Education, Okara Campus, Okara 56300, Pakistan

## Abstract
Homotopy perturbation method is used for solving the multi-point boundary value problems. The approximate solution is found in the form of a rapidly convergent series. Several numerical examples have been considered to illustrate the efficiency and implementation of the method and the results are compared with the other methods in the literature.

**K e y w o r d s**: Homotopy perturbation method; multi-point boundary value problems, linear and nonlinear Problems, approximate solution

## 1 Introduction

Multipoint boundary value problems arise in applied mathematics and physics. For example, the vibrations of a guy wire of uniform cross-section and composed of $N$ parts of different densities can be given as a multi-point boundary value problem (Moshiinsky, 1950). (Hajji, 2009), considered the multi-point boundary value problems which occurs in many areas of engineering applications such as in modelling the flow of fluid such as water, oil and gas through ground layers, where each layer constitutes a subdomain. In (Timoshenko, 1961), many problems in the theory of elastic stability are handled by multi-point problems. In (Geng and Cui, 2010) large size bridges are sometimes contrived with multi-point supports which correspond to a multi-point boundary value condition. Many authors studied the existence and multiplicity of solutions of multi-point boundary value problems (Eloe and Henderson, 2007), (Feng and Webb, 2007), (Graef and Webb, 2009), (Henderson and Kunkel, 2008), (Liu, 2003). Some research works are available on numerical analysis of the multi-point boundary value problems. Numerical solutions of multi-point boundary value problems have been studies by (Geng, 2009), (Lin and Lin, 2010), (Tatari and Dehghan, 2006), (Wu and Li, 2011). Siddiqi and Ghazala (Siddiqi and Akram, 2006a, 2006b) presented the solutions of fifth and sixth order boundary value problems using non-polynomial spline technique. Recently, Akram and Hamood (Akram and Hamood, 2013a) used the reproducing Kernel space method to solve the eighth-order boundary value problems and in (Akram and Hamood, 2013b) find the solution of a class of sixth order boundary value problems using the reproducing kernel space method. Siddiqi and Iftikhar (Siddiqi and Iftikhar, 2013) presented the solution of higher order boundary value problems using the homotopy analysis method.

J. H. He (He, 1999, 2003, 2004, 2005) developed the homotopy perturbation method for solving nonlinear initial and boundary value problems by combining the standard homotopy in topology and the perturbation technique. By this method, a rapid convergent series solution can be obtained in most of the cases. Usually, a few terms of the series solution can be used for numerical calculations. Chun C. and Sakthivel, (Chun C. and Sakthivel, 2010), implement the homotopy perturbation method for solving the linear and nonlinear two-point boundary value problems. The convergence of the homotopy perturbation method was discussed in (Biazar and H. Ghazvini, 2009), (He, 1999), (Hussein, 2011), (Turkyilmazoglu , 2011). This method has been successfully applied to ordinary differential equations, partial differential equations and other fields (Belendez, 2007), (Dehghan and Shakeri, 2008), (He, 1999, 2003, 2004, 2005), (Rana, 2007), (Yusufoglu, 2007).

_______________________________

[a] Email: shahidsiddiqiprof@yahoo.co.uk
[b] Email: miftikhar@hotmail.com, muzamil.iftikhar@ue.edu.pk

In this paper, the application of the homotopy perturbation method for finding an approximate solution for multi-point boundary value problems has been investigated.

The organization of the rest of the paper is as follows: In section 2, the homotopy perturbation method is applied to some ordinary differential equations with given multi-point boundary conditions. In section 3, the homotopy perturbation method is used to solve several examples. Finally, in section 4, the conclusion is presented.

## 2 Analysis of the Homotopy Perturbation Method(He, 1999)

Consider the nonlinear differential equation

$$L(u) + N(u) = f(r), \quad r \in \Omega \tag{2.1}$$

with boundary conditions

$$B(u, \frac{\partial u}{\partial n}) = 0, \quad r \in \Gamma \tag{2.2}$$

where $L$ is a linear operator, $N$ is a nonlinear operator, $f(r)$ is a known analytic function, $B$ is a boundary operator and $\Gamma$ is the boundary of the domain $\Omega$.

By He's homotopy perturbation technique (He, 1999), define a homotopy $v(r,p) : \Omega \times [0,1] \to R$ which satisfies

$$H(v, p) = (1-p)[L(v) - L(u_0)] + p[L(v) + N(v) - f(r)] = 0, \tag{2.3}$$

or

$$H(v, p) = L(v) - L(u_0) + pL(u_0) + p[N(v) - f(r)] = 0, \tag{2.4}$$

where $r \in \Omega$, $p \in [0,1]$ is an embedding parameter and $u_0$ is an initial approximation of Eq. ( 2. 1 ) which satisfies the boundary conditions. Clearly

$$H(v, 0) = L(v) - L(u_0) = 0, \tag{2.5}$$

$$H(v, 1) = L(v) + N(v) - f(r) = 0, \tag{2.6}$$

As $p$ changes from $0$ to $1$, then $v(r, p)$ changes from $u_0(r)$ to $u(r)$ This is called a deformation and $L(v) - L(u_0)$, $L(v) + N(v) - f(r)$ are said to be homotopic in topology. According to the homotopy perturbation method, firstly, the embedding parameter $p$ can be used as a small parameter, and assume that the solution of Eq. ( 2. 3 ) and Eq. ( 2. 4 ) can be expressed as a power series in $p$, that is,

$$v = v_0 + pv_1 + p^2 v_2 + \cdots \tag{2.7}$$

For $p = 1$, the approximate solution of Eq. ( 2. 1 ) therefore, can be expressed as

$$v = \lim_{p \to 1} v = v_0 + v_1 + v_2 + \cdots \tag{2.8}$$

The series in Eq. ( 2. 8 ) is convergent in most cases, and the convergence rate of the series depends on the nonlinear operator, see (Biazar and H. Ghazvini, 2009), (He, 1999). Moreover, the following judgments are made by He ( He, 1999), ( He, 2006),

(i) The second order derivative of $N(v)$ w.r.t. $v$ must be small as the parameter may be reasonably large, i.e., $p \to 1$.

(ii) $\left\| L^{-1}\left(\frac{\partial N}{\partial v}\right) \right\|$ must be smaller than one, so that, the series converges.

To implement the method, several numerical examples are considered in the following section.

## 3 Numerical Examples

**Example 3.1** Consider the following third-order linear differential equation with three point boundary conditions

$$\left. \begin{array}{l} u'''(x) - k^2 u'(x) + a = 0, 0 \leq x \leq 1, \\ u'(0) = u'(1) = 0, u(0.5) = 0. \end{array} \right\} \quad (3.1)$$

The exact solution of the Example 3.1 is

$$u(x) = \frac{a}{k^3}(\sinh\frac{k}{2} - \sinh kx) + \frac{a}{k^2}(x - \frac{1}{2}) + \frac{a}{k^3}\tanh\frac{k}{2}(\cosh kx - \cosh\frac{k}{2})$$

where the constants are $k = 5$ and $a = 1$ (Akram et al., 2013), (Ali et al., 2010), (Saadatmandi and Dehghan, 2012), (Tirmizi et al., 2005).

Using the homotopy perturbation method, the following homotopy for the system ( 3. 1 ) is constructed

$$u''' = p[25u'] - 1, \quad (3.2)$$

where $p \in [0,1]$ is the embedding parameter. Assume that the solution of Problem ( 3. 1 ) is

$$u = u_0 + pu_1 + p^2 u_2 + \cdots \quad (3.3)$$

Substituting Eq. ( 3. 3 ) in Eq. ( 3. 2 ), and equating the coefficients of like powers of $p$, gives the following set of differential equations

$$p^0: \quad u_0''' = -1, \ u_0'(0) = 0, u_0(0) = A, u_0''(0) = B$$
$$p^1: \quad u_1''' = 25u_0', \ u_1'(0) = 0, u_1(0) = 0, u_1''(0) = 0$$
$$p^2: \quad u_2''' = 25u_1', \ u_2'(0) = 0, u_2(0) = 0, u_2''(0) = 0$$
$$\vdots$$

where $A$ and $B$ are unknown constants to be determined. The corresponding solutions for the above system of equations are the series solution given as

$$u_0(x) = \frac{1}{6}(6A + 3Bx^2 - x^3)$$
$$u_1(x) = \frac{5}{24}(5Bx^4 - x^5)$$
$$\vdots$$

Using the 11-term approximation, that is

$$U(x) = u_0(x) + u_1(x) + u_2(x) + \cdots + u_{10}(x) \quad (3.4)$$

Imposing the boundary conditions of the system ( 3. 1 ) on Eq. ( 3. 4 ) the values of the constants A and B can be obtained as

$A = -0.012107085822126442,$ $\qquad B = 0.19732286064025403.$

Then, the series solution can be expressed as

$U(x) = -0.0121071 + 0.0986614x^2 - 0.16667x^3 + 0.205545x^4 - 0.208333x^5 + 0.171287x^6$
$\quad - 0.124008x^7 + 0.0764675x^8 - 0.0430583x^9 + 0.021241x^{10} - 0.009785x^{11}$
$\quad + 0.00402291x^{12} + O(x^{13}).$

The comparison of the approximate series solution of the problem (3.1) with the results of methods in (Akram et al., 2013), (Ali et al., 2010), (Saadatmandi and Dehghan, 2012), (Tirmizi et al., 2005) is given in Table 1, which shows that the method is quite efficient. In Figure 1(a) and Figure 1(b) errors $|U - u_{Exact}|$ and $|U - u_{Exact}|/|u_{Exact}|$ are plotted respectively. Figure 1 shows that the method is in excellent agreement with (Tatari and Dehghan, 2007).

**Example 3.2** Consider the linear fourth-order nonlocal boundary value problem

$$\left.\begin{array}{l} u^{(4)}(x) - e^x u^{(3)}(x) + u(x) = 1 - e^x \cosh(x) + 2\sinh(x), 0 \le x \le 1 \\ u\left(\dfrac{1}{4}\right) = 1 + \sinh\left(\dfrac{1}{4}\right), u^{(1)}\left(\dfrac{1}{4}\right) = 1 + \cosh\left(\dfrac{1}{4}\right), \\ u^{(2)}\left(\dfrac{1}{4}\right) = \sinh\left(\dfrac{1}{4}\right), u\left(\dfrac{1}{2}\right) - u\left(\dfrac{3}{4}\right) = \sinh\left(\dfrac{1}{2}\right) - \sinh\left(\dfrac{3}{4}\right). \end{array}\right\}$$

The exact solution of the problem (3.2) is $u(x) = 1 + \sinh(x)$ (Lin and Lin, 2010), (Wu and Li, 2011).

Using the homotopy perturbation method, the following homotopy for the system ( 3. 5 ) is constructed

$$u^{(4)}(x) = 1 - e^x \cosh(x) + 2\sinh(x) + p[e^x u^{(3)}(x) - u(x)], \qquad (3.6)$$

where $p \in [0,1]$ is the embedding parameter. Assume that the solution of Problem ( 3. 5 ) is

$$u = u_0 + p u_1 + p^2 u_2 + \cdots \qquad (3.7)$$

Substituting Eq. ( 3. 7 ) in Eq. ( 3. 6 ), and equating the coefficients of like powers of $p$, gives the following set of differential equations

$p^0:$ $\quad u_0^{(4)}(x) = 1 - e^x \cosh(x) + 2\sinh(x), u_0(0) = A, u_0^{(1)}(0) = B, u_0^{(2)}(0) = C, u_0^{(3)}(0) = D,$

$p^1:$ $\quad u_1^{(4)}(x) = e^x u_0^{(3)} - u_0, u_1(0) = 0, u_1^{(1)}(0) = 0, u_1^{(2)}(0) = 0, u_1^{(3)}(0) = 0,$

$p^2:$ $\quad u_2^{(4)}(x) = e^x u_1^{(3)} - u_1, u_2(0) = 0, u_2^{(1)}(0) = 0, u_2^{(2)}(0) = 0, u_2^{(3)}(0) = 0,$

$\vdots$

where $A$, $B$, $C$ and $D$ are unknown constants to be determined. The corresponding solutions for the above system of equations are the series solution given as

$u_0(x) = \dfrac{1}{96}(-96 + 96e^{2x} - 3e^{3x} + e^x(3 + 96A + 6(-31 + 16B)x + 6(1 + 8C)x^2 + 4(-7 + 4D)x^3 + 2x^4))$

$u_1(x) = \dfrac{1}{1451520}(e^{-x}(1451520 + 93555e^{3x} - 4480e^{4x} + 362880e^{2x}(-19 + 4C + 2x) + e^x(5354125$

$+7446810x+1828890x^2+923580x^3-1890(-31+32A)x^4-756(-31+16B)x^5-252(1+8C)x^6+504x^7-18x^8-288D(5040+5040x+2520x^2+840x^3+x^7))))$

$\vdots$

Using only 6-term approximation, that is

$$U(x)=u_0(x)+u_1(x)+u_2(x)+\cdots+u_5(x) \quad (3.8)$$

Imposing the boundary conditions of the system (3.5) on Eq. (3.8) the values of the constants $A$, $B$, $C$ and $D$ can be obtained as

$A = 0.9999999980259633$, $\qquad B = 1.0000000216759806$,

$C = -1.6366491839105507 \times 10^{-7}$, $\qquad D = 1.00000056811826$.

Then, the series solution can be expressed as

$U(x) = 1 + x - 8.16726 \times 10^{-8} x^2 + 0.16667 x^3 + 2.38316 \times 10^{-8} x^4 + 0.00833334 x^5 + 4.15407 \times 10^{-9} x^6$

$+0.000198414 x^7 + 7.15518 \times 10^{-10} x^8 \, 2.75604 \times 10^{-6} x^9 + 1.2948 \times 10^{-10} x^{10} - 02.49988$

$\times 10^{-8} x^{11} - 1.31503 \times 10^{-8} x^{12} + O(x^{13}).$ (3.9)

The approximate series solution of the problem (3.2) is compared with $u(x) = 1 + \sinh(x)$ (Lin and Lin, 2010), (Wu and Li, 2011) in Table 2, which shows that the method is quite efficient. Absolute errors $|U - u_{Exact}|$ are plotted in Figure 2.

**Example 3.3** The following fourth order nonlinear boundary value problem is considered

$$\left.\begin{array}{l} u^{(4)}(x) - e^{-x} u^2(x) = 0, 0 \leq x \leq 1 \\ u(0) = u^{(1)}(0) = 1, u\left(\dfrac{3}{4}\right) = e^{\frac{3}{4}}, u(1) = e. \end{array}\right\} \quad (3.10)$$

The exact solution of the problem (3.3) is $u(x) = e^x$.

Using the homotopy perturbation method, the following homotopy for the system (3.10) is constructed

$$u^{(4)}(x) = p[e^{-x} u^2], \quad (3.11)$$

where $p \in [0,1]$ is the embedding parameter. Assume that the solution of the given problem is

$$u = u_0 + p u_1 + p^2 u_2 + \cdots \quad (3.12)$$

The nonlinear term $N(u)$ in Eq. (3.11) can be expressed as

$$N(u) = N(u_0) + p N(u_0, u_1) + p^2 N(u_0, u_1, u_2) + \cdots, \quad (3.13)$$

Where

$$N(u_0, u_1, \ldots, u_n) = \dfrac{1}{n!} \dfrac{d^n}{dp^n}\left[N\left(\sum_{k=0}^{n} p^k u_k\right)\right]_{p=0}, \quad n = 0, 1, 2, \ldots$$

is called He's polynomial (Ghorbani, 2009) Substituting Eq. (3.12) and Eq. (3.13) in Eq. (3.11), and equating the coefficients of like powers of $p$, gives the following set of differential equations

$p^0: \quad u_0^{(4)}(x) = 0, \ u_0(0) = 1, u_0^{(1)}(0) = 1, u_0^{(2)}(0) = A, u_0^{(3)}(0) = B,$

$p^1: \quad u_1^{(4)}(x) = e^{-x} u_0^2, \ u_1(0) = 0, u_1^{(1)}(0) = 0, u_1^{(2)}(0) = 0, u_1^{(3)}(0) = 0,$

$p^2: \quad u_2^{(4)}(x) = e^{-x}u_0u_1, \ u_2(0) = 0, u_2^{(1)}(0) = 0, u_2^{(2)}(0) = 0, u_2^{(3)}(0) = 0,$
$\vdots$

where $A$ and $B$ are unknown constants to be determined. Following Example (3.1), using the 3-term approximation and imposing the boundary conditions at $x = 0.75$ and $x = 1$, the constants are obtained as

$A = 0.9999994087690695, \qquad B = 1.0000024198861392.$

Then, the series solution can be expressed as

$U(x) = 1 + x - 0.5x^2 + 0.166667x^3 + 0.416667x^4 + 0.008333337x^5 + 0.00138889x^6$
$\quad + 0.000198414x^7 + 0.0000248016x^8 + 2.75573 \times 10^{-6} x^9 + 0.75571 \times 10^{-7} x^{10}$
$\quad + 2.50527 \times 10^{-8} x^{11} - 1.524 \times 10^{-7} x^{12} + O(x^{13}).$

In Table 3, the comparison of the exact solution with the series solution of the problem (3.3) is given, which shows that the method is quite efficient. In Figure 3 absolute errors $|U - u_{Exact}|$ are plotted in Figure 3.

**Example 3.4** The following fifth order nonlinear three points boundary value problem is considered

$$u^{(5)}(x) - e^{-x}u^2(x) = 0, 0 < x < 1$$
$$u(0) = u^{(1)}(0) = 1, u\left(\frac{1}{2}\right) = e^{\frac{1}{2}}, u(1) = u^{(2)}(1) = e. \qquad (3.14)$$

The exact solution of the problem (3.4) is $u(x) = e^x$.

Using the homotopy perturbation method, the following homotopy for the system (3.14) is constructed

$$u^{(5)}(x) = p[e^{-x}u^2], \qquad (3.15)$$

where $p \in [0,1]$ is the embedding parameter. Assume that the solution of the given problem is

$$u = u_0 + pu_1 + p^2 u_2 + \cdots \qquad (3.16)$$

The nonlinear term $N(u)$ in Eq. (3.11) can be expressed as

$$N(u) = N(u_0) + pN(u_0, u_1) + p^2 N(u_0, u_1, u_2) + \cdots, \qquad (3.17)$$

Where

$$N(u_0, u_1, \ldots, u_n) = \frac{1}{n!}\frac{d^n}{dp^n}\left[N\left(\sum_{k=0}^{n} p^k u_k\right)\right]_{p=0}, \quad n = 0, 1, 2, \ldots$$

is called He's polynomial (Ghorbani, 2009). Substituting Eq. (3.16) and Eq. (3.17) in Eq. (3.15), and equating the coefficients of like powers of $p$, gives the following set of differential equations

$p^0: \quad u_0^{(5)}(x) = 0, \ u_0(0) = 1, u_0^{(1)}(0) = 1, u_0^{(2)}(0) = A, u_0^{(3)}(0) = B, u_0^{(4)}(0) = C,$
$p^1: \quad u_1^{(5)}(x) = e^{-x}u_0^2, \ u_1(0) = 0, u_1^{(1)}(0) = 0, u_1^{(2)}(0) = 0, u_1^{(3)}(0) = 0, u_1^{(4)}(0) = 0,$
$p^2: \quad u_2^{(5)}(x) = e^{-x}u_0u_1, \ u_2(0) = 0, u_2^{(1)}(0) = 0, u_2^{(2)}(0) = 0, u_2^{(3)}(0) = 0, u_2^{(4)}(0) = 0,$
$\vdots$

where $A$, $B$ and $C$ are unknown constants to be determined. Following Example (3.1), using the 3-term approximation and imposing the boundary conditions at $x = 0.75$ and $x = 1$, the constants are obtained as

$$A = 1.00000000568, \quad B = 0.99999994805, \quad C = 1.00000014256.$$

Then, the series solution can be expressed as

$$U(x) = 1 + x - 0.5000000028x^2 + 0.166666x^3 + 0.4166667x^4 + 0.008333333x^5 + 0.001388889x^6$$
$$+ 0.00019841x^7 + 0.00002480x^8 + 2.7557327 \times 10^{-6} x^9 + 2.7557319 \times 10^{-7} x^{10}$$
$$+ 2.50521 \times 10^{-8} x^{11} + 2.087675 \times 10^{-9} x^{12} + O(x^{13}).$$

In Table 4, the comparison of the exact solution with the series solution of the problem (3.4) is given, which shows that the method is quite efficient. In Figure 4 absolute errors $|U - u_{Exact}|$ are plotted.

**Example 3.5** The following sixth order nonlinear boundary value problem is considered

$$\left. \begin{array}{l} u^{(6)}(x) - e^{-x} u^2(x) = 0, 0 < x < 1 \\ u(0) = u^{(1)}(0) = u^{(2)}(0) = u^{(3)}(0) = 1, u\left(\dfrac{1}{2}\right) = e^{\frac{1}{2}}, u(1) = e. \end{array} \right\} \quad (3.18)$$

The exact solution of the problem (3.5) is $u(x) = e^x$.

Using the aforesaid method, the series solution can be expressed as

$$U(x) = 1 + (1.)x - 0.5x^2 + 0.166667x^3 + 0.4166667x^4 + 0.008333333x^5 + 0.00138885x^6$$
$$+ 0.000198432x^7 + 0.0000247952x^8 + 2.75728 \times 10^{-6} x^9 + 2.75381 \times 10^{-7} x^{10}$$
$$+ 2.49973 \times 10^{-8} x^{11} + 2.14303 \times 10^{-9} x^{12} + O(x^{13}).$$

The comparison of the exact solution with the series solution of the problem (3.5) is given in Table 5, which shows that the method is quite accurate.

**Example 3.6** The following seventh order nonlinear boundary value problem is considered

$$\left. \begin{array}{l} u^{(7)}(x) - e^{-x} u^2(x) = 0, 0 < x < 1 \\ u(0) = u^{(1)}(0) = u^{(2)}(0) = u^{(3)}(0) = u^{(4)}(0) = 1, u\left(\dfrac{1}{2}\right) = e^{\frac{1}{2}}, u(1) = e. \end{array} \right\} \quad (3.19)$$

The exact solution of the problem (3.6) is $u(x) = e^x$.

Using the aforesaid method, the series solution can be expressed as

$$U(x) = 0.999998 + (1.)x - 0.499998x^2 + 0.166668x^3 + 0.4166661x^4 + 0.00833367x^5 + 0.00138876x^6$$
$$+ 0.000198417x^7 + 0.0000248361x^8 + 2.72677 \times 10^{-6} x^9 + 2.89152 \times 10^{-7} x^{10}$$
$$+ 2.14384 \times 10^{-8} x^{11} + 2.0249 \times 10^{-9} x^{12} + O(x^{13}).$$

The comparison of the exact solution with the series solution of the problem (3.6) is given in Table 6, which shows that the method is quite accurate.

**Example 3.7** The following seventh order nonlinear boundary value problem is considered

$$u^{(7)}(x) = -u(x) - e^x(35 + 12x + 2x^2), 0 \leq x \leq 1$$
$$u(0) = 0, u^{(1)}(0) = 1, u^{(2)}(0) = 0, u^{(3)}(0) = -3, u^{(4)}(0) = -8, \quad (3.20)$$
$$u\left(\frac{1}{2}\right) = \frac{e^{\frac{1}{2}}}{4}, u(1) = e.$$

The exact solution of the problem (3.7) is $u(x) = x(1-x)e^x$.
Using the aforesaid method, the series solution can be expressed as
$$U(x) = x - 0.5x^3 - 0.333333x^4 - 0.125x^5 - 0.333333x^6 - 0.00694444x^7 - 0.00119048x^8$$
$$-0.000173611x^9 - 0.0000220459x^{10} - 2.48016 \times 10^{-6} x^{11} - 2.50521 \times 10^{-7} x^{12} + O(x^{13}).$$

The comparison of the exact solution with the series solution of the problem (3.7) is given in Table 7, which shows that the method is quite accurate.

**Conclusion** In this paper, the homotopy perturbation method has been applied to solve the multi-point boundary value problems. It is clearly seen that homotopy method is a powerful and accurate method for finding solutions for multi-point boundary value problems in the form of analytical expressions and presents a rapid convergence for the solutions. The numerical results showed that the homotopy perturbation method can solve the problem effectively and the comparison shows that the present method is in good agreement with the existing results in the literature.

Table 1: Comparison of numerical results for Example 3.1

| $x$ | Exact solution | Approximate Series solution | Absolute Error Present method | Absolute Error (Tirmizi et al., 2005) | Absolute Error (Ali et al., 2010) | Absolute Error (Akram et al., 2013) |
|---|---|---|---|---|---|---|
| 0.0 | -0.0121071 | -0.0121071 | 2.07338E-10 | 0.00003515 | 1.298 E-10 | 8.37E-07 |
| 0.1 | -0.0112685 | -0.0112685 | 2.02182E-10 | 0.00003850 | 3.099E-09 | 3.39E-07 |
| 0.2 | -0.00922221 | -0.00922221 | 1.85398E-10 | 0.00003028 | 6.959E-09 | 9.16E-08 |
| 0.3 | -0.00646687 | -0.00646687 | 1.52702E-10 | 0.00002231 | 1.086E-09 | 7.22E-08 |
| 0.4 | -0.00332019 | -0.00332019 | 9.57487E-11 | 0.00001403 | 1.065E-08 | 7.86E-08 |
| 0.5 | 0.000000 | -4.03581E-18 | 4.03581E-18 | 0.00000700 | 6.155E-17 | 6.55E-08 |
| 0.6 | 0.00332019 | 0.00332019 | 1.58981E-10 | 0.00001260 | 1.065E-08 | 6.35E-08 |
| 0.7 | 0.00646687 | 0.00646687 | 4.21657E-10 | 0.00001260 | 1.086E-09 | 6.26E-08 |
| 0.8 | 0.00922221 | 0.00922221 | 8.52404E-10 | 0.00001956 | 6.959E-09 | 9.54E-08 |
| 0.9 | 0.0112685 | 0.0112685 | 1.51972E-09 | 0.00002741 | 3.099E-09 | 3.37E-07 |
| 1.0 | 0.0121071 | 0.0121071 | 2.1212E-09 | 0.00002395 | 1.298E-10 | 8.48E-07 |

Table 2: Comparison of numerical results for problem (3.2)

| $x$ | Exact solution | Approximate series solution | Absolute Error present method | Absolute Error in (Lin and Lin, 2010) | Absolute Error in (Wu and Li, 2011) |
|---|---|---|---|---|---|
| 0.0 | 1.0000 | 1.0000 | 1.95677E-09 | 1.02E-4 | 2.54E-8 |
| 0.1 | 1.10017 | 1.10017 | 5.83738E-10 | 1.81E-5 | 4.70E-9 |
| 0.2 | 1.20134 | 1.20134 | 1.04897E-10 | 5.33E-7 | 1.39E-10 |
| 0.3 | 1.30452 | 1.30452 | 9.55858E-11 | 3.94E-7 | 1.25E-10 |
| 0.4 | 1.41075 | 1.41075 | 3.01129E-10 | 7.60E-6 | 2.40E-9 |
| 0.5 | 1.5211 | 1.5211 | 2.05711E-09 | 2.36E-5 | 7.58E-9 |
| 0.6 | 1.63665 | 1.63665 | 5.14444E-09 | 3.90E-5 | 1.13E-8 |
| 0.7 | 1.75858 | 1.75858 | 6.73259E-09 | 3.73E-5 | 4.30E-9 |
| 0.8 | 1.88811 | 1.88811 | 1.47844E-08 | 2.42E-6 | 2.80E-8 |
| 0.9 | 2.02652 | 2.02652 | 1.55269E-07 | 1.06E-4 | 1.05E-7 |
| 1.0 | 2.1752 | 2.1752 | 7.92993E-07 | 3.05E-4 | 2.52E-7 |

Table 3: Comparison of numerical results for problem (3.3)

| $x$ | Exact solution | Approximate series solution | Absolute Error present method |
|---|---|---|---|
| 0.0 | 1.00000 | 1.0000 | 6.26543E-12 |
| 0.1 | 1.10517 | 1.10517 | 2.55342E-09 |
| 0.2 | 1.2214 | 1.2214 | 8.60246E-09 |
| 0.3 | 1.34986 | 1.34986 | 1.57141E-08 |
| 0.4 | 1.49182 | 1.49182 | 2.1502E-08 |
| 0.5 | 1.64872 | 1.64872 | 2.35332E-08 |
| 0.6 | 1.82212 | 1.82212 | 1.96291E-08 |
| 0.7 | 2.01375 | 2.01375 | 8.27396E-09 |
| 0.8 | 2.22554 | 2.22554 | 9.18081E-09 |
| 0.9 | 2.4596 | 2.4596 | 2.28539E-08 |
| 1.0 | 2.71828 | 2.71828 | 8.86402E-12 |

Table 4: Comparison of numerical results for problem (3.4)

| $x$ | Exact solution | Approximate series solution | Absolute Error present method |
|---|---|---|---|
| 0.0 | 1.00000 | 1.00000 | 0.0000 |
| 0.1 | 1.10517 | 1.10517 | 5.58569E-10 |
| 0.2 | 1.2214 | 1.2214 | 3.80139E-10 |
| 0.3 | 1.34986 | 1.34986 | 4.5143E-10 |
| 0.4 | 1.49182 | 1.49182 | 2.60672E-10 |
| 0.5 | 1.64872 | 1.64872 | 2.39371E-10 |
| 0.6 | 1.82212 | 1.82212 | 7.77565E-11 |
| 0.7 | 2.01375 | 2.01375 | 1.64396E-10 |
| 0.8 | 2.22554 | 2.22554 | 8.80967E-10 |
| 0.9 | 2.4596 | 2.4596 | 1.3927E-10 |
| 1.0 | 2.71828 | 2.71828 | 2.4848E-10 |

Table 5: Comparison of numerical results for problem (3.5)

| x | Exact solution | Approximate series solution | Absolute Error present method |
|---|---|---|---|
| 0.0 | 1.00000 | 1.00000 | 7.77951E-09 |
| 0.1 | 1.10517 | 1.10517 | 1.16784E-08 |
| 0.2 | 1.2214 | 1.2214 | 7.57914E-09 |
| 0.3 | 1.34986 | 1.34986 | 2.04205E-08 |
| 0.4 | 1.49182 | 1.49182 | 1.75262E-08 |
| 0.5 | 1.64872 | 1.64872 | 1.03601E-08 |
| 0.6 | 1.82212 | 1.82212 | 1.60579E-09 |
| 0.7 | 2.01375 | 2.01375 | 4.20526E-10 |
| 0.8 | 2.22554 | 2.22554 | 2.25408E-08 |
| 0.9 | 2.4596 | 2.4596 | 8.26443E-09 |
| 1.0 | 2.71828 | 2.71828 | 1.69864E-08 |

Table 6: Comparison of numerical results for problem (3.6)

| x | Exact solution | Approximate series solution | Absolute Error present method |
|---|---|---|---|
| 0.0 | 1.00000 | 1.00000 | 8.57339E-07 |
| 0.1 | 1.10517 | 1.10517 | 1.02848E-06 |
| 0.2 | 1.2214 | 1.2214 | 4.81648E-07 |
| 0.3 | 1.34986 | 1.34986 | 4.94121E-06 |
| 0.4 | 1.49182 | 1.49182 | 5.39926E-07 |
| 0.5 | 1.64872 | 1.64872 | 6.97339E-07 |
| 0.6 | 1.82212 | 1.82212 | 5.62541E-07 |
| 0.7 | 2.01375 | 2.01375 | 3.98414E-07 |
| 0.8 | 2.22554 | 2.22554 | 4.86015E-07 |
| 0.9 | 2.4596 | 2.4596 | 1.5371E-06 |
| 1.0 | 2.71828 | 2.71828 | 6.25725E-07 |

Table 7: Comparison of numerical results for Example 3.7

| x | Exact solution | Approximate series solution | Absolute Error Present method |
|---|---|---|---|
| 0.0 | 0.0000 | 0.0000 | 0.0000 |
| 0.1 | 0.9946 | 0.9946 | 5.69961E-14 |
| 0.2 | 0.1954 | 0.1954 | 8.9373E-15 |
| 0.3 | 0.2835 | 0.2835 | 4.05231E-15 |
| 0.4 | 0.3580 | 0.3580 | 1.54876E-14 |
| 0.5 | 0.4122 | 0.4122 | 1.4555E-13 |
| 0.6 | 0.4373 | 0.4373 | 1.03195E-13 |
| 0.7 | 0.4229 | 0.4229 | 4.16889E-14 |
| 0.8 | 0.3561 | 0.3561 | 2.33036E-13 |
| 0.9 | 0.2214 | 0.2214 | 2.39697E-13 |
| 1.0 | 0.0000 | -2.12607E-13 | 2.12607E-13 |

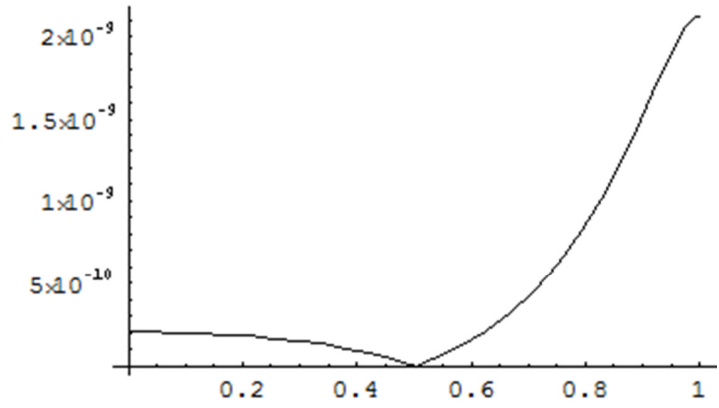

Figure 1(a): Plot of errors $|U - u_{Exact}|$.

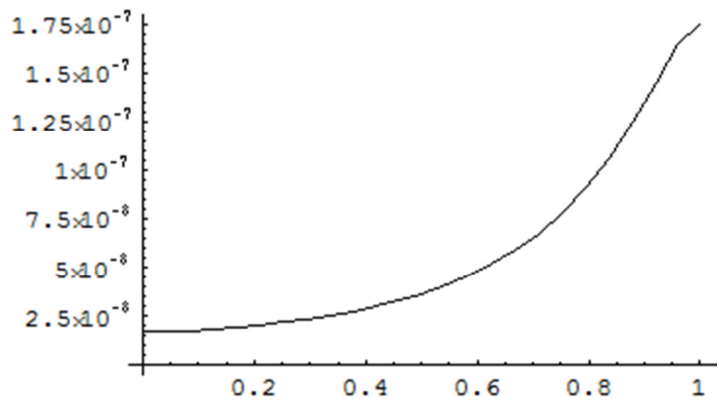

Figure 1(b): Plot of errors $|U - u_{Exact}|/|u_{Exact}|$.

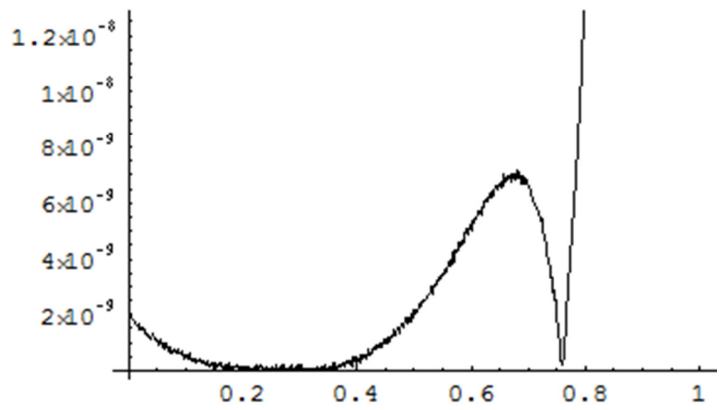

Figure 2: Plot of errors $|U - u_{Exact}|$.

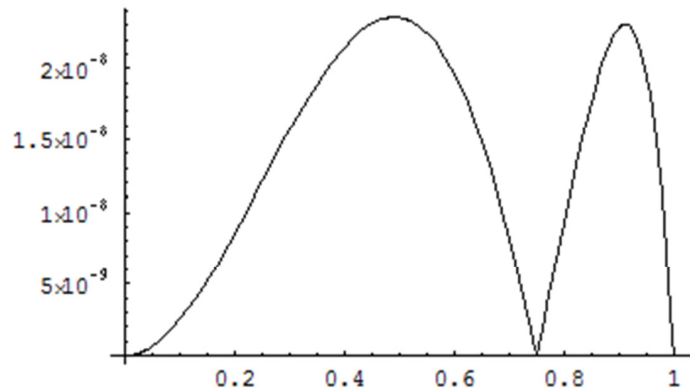

Figure 3: Plot of absolute errors $|U - u_{Exact}|$.

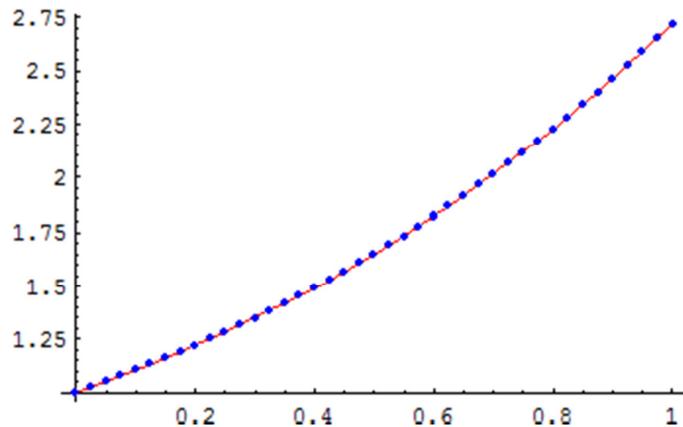

Figure 4: Comparison of the approximate solution with the exact solution for problem (3.4). Dotted line: approximate solution, solid line: the exact solution.